\documentclass[10pt]{article}

\usepackage{amsfonts,amsmath,amscd, amssymb,  vmargin}
\usepackage[francais]{babel}
\usepackage[all]{xy} 
\everymath{\displaystyle}

\begin{document} 
\title{Un scindage de l'application de Frobenius sur toute  l'alg\`ebre  des distributions    de $SL_{2}$      }
\author{Michel Gros  
\\IRMAR, UMR CNRS 6625\\Universit\'e de Rennes I\\Campus de Beaulieu\\35042 Rennes
cedex\\France\\e-mail : michel.gros@univ-rennes1.fr}

\maketitle
\author

{\it{Abstract}}. We define, over ${\Bbb{F}}_{p}$  ($p>2$),  a splitting of the Frobenius morphism $Fr :  {\text{Dist}}\,(G) \rightarrow {\text{Dist}}\,(G)$ on the whole ${\text{Dist}}\,(G)$, the algebra of distributions of the $k$-algebraic group   $G:=SL_{2}$. This splitting is compatible (and lifts) the theory of Frobenius descent for arithmetic ${\cal{D}}$-modules over $X:={\Bbb{P}}_{k}^{1}$.

\setcounter{section}{-1}
\section{Introduction} 

Soient $p > 2$ un nombre premier, $k = {\Bbb{F}}_{p}$ le corps fini \`a $p$ \'el\'ements,  $X$ un $k$-sch\'ema lisse et $({\cal{D}}_{X}^{(m)}) _{m \in {\Bbb{N}}}$ le syst\`eme inductif d' anneaux d' op\'erateurs diff\'erentiels d\'efini par Berthelot (\cite {Berthelot1}, 2. ou \cite{Berthelot3} pour un survol de toute la th\'eorie). Si $F  : X \rightarrow X'$ d\'esigne le morphisme de Frobenius relatif attach\'e \`a cette situation, un r\'esultat fondamental (\cite {Berthelot2}, 2.3.6 et 2.4.6) est que le foncteur  $F ^{*} : ({\cal{D}}_{X'}^{(m)}-{\text{Mod}}) \rightarrow ({\cal{D}}_{X}^{(m+1)}-{\text{Mod}})$ est  une \'equivalence de cat\'egories. La question, non triviale,  de l'explicitation locale (\og descente par Frobenius\fg ), partant d' un  ${\cal{D}}_{X}^{(m+1)}$-module du ${\cal{D}}_{X'}^{(m)}$-module qui lui correspond a \'et\'e \'etudi\'ee dans  \cite{Garnier} et   un certain op\'erateur diff\'erentiel "projecteur"\footnote{Nous nous \'ecartons ici de la typographie  de ${\it{loc. cit.}}$  dans laquelle ce projecteur est not\'e $H$ afin d'\'eviter plus bas toute confusion avec  la notation standard de la base canonique de $sl_2$.}  ${\cal{H}}$ (cf. ${\it{loc.\, cit.}}$  2.5)  joue un r\^ole central.\\

Si, d'autre part, $G$ d\'esigne  un $k$-groupe r\'eductif   d'alg\`ebre de Lie ${\cal{G}}$, l'id\'ee, fondamentale, d'introduction  des "puissances divis\'ees partielles"  dans  la d\'efinition des  $({\cal{D}}_{X}^{(m)}) _{m \in {\Bbb{N}}}$ une fois appliqu\'ee \`a $G$  permet de d\'efinir  (cf. \cite{Kaneda-Ye}) un syst\`eme inductif d'alg\`ebres enveloppantes $(U^{(m)}({\cal{G}})  ) _{m \in {\Bbb{N}}}$ tel que $U^{(0)}({\cal{G}})  $ s'identifie \`a  l'alg\`ebre enveloppante universelle et la limite inductive des $U^{(m)}({\cal{G}})  $ s'identifie \`a l'alg\`ebre des distributions ${\text{Dist}}\,(G)$ de $G$. Notant $X (={\Bbb{P}}^{1}_{k})$ le $k$-sch\'ema des drapeaux (complets) de $G$, on dispose de morphismes\footnote{Variante "\`a niveau" de \cite{Kashiwara}, Prop. 6.2.3.}  d'anneaux $\rho_{m} : U^{(m)}({\cal{G}})    \rightarrow H^{0}(X,  {\cal{D}}_{X}^{(m)})$ (et d'une th\'eorie de la localisation\footnote{Le r\'esultat le plus complet concerne le cas $m=0$ (\cite{Bezrukavnikov}) et des r\'esultats partiels existent pour $m > 0$ (\cite{Kaneda}).} permettant  $\it{in\, fine}$ de passer des $U^{(m)}({\cal{G}})  $-modules aux  ${\cal{D}}_{X}^{(m)}$-modules).\\

Nous abordons dans ce travail,  dans le cas $G =  SL_{2} $, la question suivante :    existe-t-il    des variantes de la   "descente par Frobenius" pour les $U^{(m)}({\cal{G}})  $-modules  (ou pour certains d' entre eux)  compatibles (via $\rho_{m}$) avec ce qui existe pour les ${\cal{D}}_{X}^{(m)}$-modules ?  Le r\'esultat principal\footnote{Nous ne donnons la d\'emonstration compl\`ete que d' un r\'esultat un peu moins pr\'ecis.}   est le suivant :\\

{\bf{Proposition}} (cf. Prop. 2.2.1). Il existe, pour tout $m \geq 0$ un morphisme (non unif\`ere) d'alg\`ebres $\varphi_{m} : U^{(m)}({\cal{G}}) \rightarrow U^{(m+1)}({\cal{G}})$ scindant le morphisme (canonique) de Frobenius $Fr : U^{(m+1)}({\cal{G}}) \rightarrow U^{(m)}({\cal{G}})$.\\

On trouve dans la litt\'erature un scindage  de $Fr$ d\'efini (dans cadre plus g\'en\'eral) par Lusztig  (et utilis\'e dans  \cite {Kumar2}, 2., \cite {Kumar1},...)   seulement sur   ${\text{Dist}}\,(B) \subset {\text{Dist}}\,(G)$ avec $B$   un Borel de $G$ et qui ne s'\'etend pas \`a  ${\text{Dist}}\,(G)$. Le caract\`ere un peu surprenant  de ce r\'esultat est que   la question de l'existence d' un rel\`evement (via $\rho_{m}$) de ${\cal{H}}$ au niveau des alg\`ebres enveloppantes, harmonisant en quelque sorte    \cite{Haboush} et     \cite{Garnier},  conduit naturellement \`a un scindage de $Fr$ sur ${\text{Dist}}\,(G)$ tout entier.\\   

Dans la partie I, nous avons choisi, afin d'\^etre le plus bref possible,  de pr\'esenter le contexte arithm\'etique des op\'erateurs diff\'erentiels et des alg\`ebres enveloppantes sous une forme ${\it{ad\, hoc}}$ (ie. par "g\'en\'erateurs et relations"). Dans la partie II, apr\`es avoir \'etabli  quelques propri\'et\'es suppl\'ementaires de la "norme" introduite dans \cite{Haboush}, nous d\'efinissons le scindage \'evoqu\'e dans la proposition ci-dessus. Que ce soit un morphisme d' alg\`ebres est \'etabli par une v\'erification reposant sur une combinaison  de  congruences \'el\'ementaires et sur les propri\'et\'es de la "norme".  Enfin, dans la partie III, nous examinons la compatibilit\'e de ce scindage avec l'application "canonique" de Berthelot et le projecteur ${\cal{H}}$ de Garnier.\\

\section{Rappels} 

\subsection{Op\'erateurs diff\'erentiels de niveau $m$} 

Soit $X$ un $k$-sch\'ema lisse. On renvoie \`a \cite{Berthelot1} et \cite{Berthelot2} pour un expos\'e complet des fondements "naturels" (faisceau des parties principales de niveau $m$,...) de la th\'eorie des $ {\cal{D}}_{X}^{(m)}$-modules   (avec ${m \in {\Bbb{N}}}$). Pour ce dont nous aurons besoin, il suffit d'adopter la pr\'esentation suivante (cf. \cite{Kaneda}). On filtre l' anneau des op\'erateurs diff\'erentiels "usuels" (ie. ceux de EGA IV) par l' ordre : ${\cal{D}}_{X} = \bigcup_{n\in {\Bbb{N}}} {\cal{D}}_{X}^{n}$ et l'on pose pour $ m \geq 0$ :\\

${\cal{D}}_{X}^{(m)} : = {\Bbb{T}}_{k} (  {\cal{D}}_{X}^{2p^{m}-1}) / (\lambda-\lambda1_{O_{X}}, (D\otimes D'- D'\otimes D) -[D,D'] , D\otimes D''-DD'', \lambda \in k, D'' \in  {\cal{D}}_{X}^{ p^{m}-1}; D,D' \in  {\cal{D}}_{X}^{ p^{m} })$\\

avec ${\Bbb{T}}_{k}(.)$ d\'esignant l'alg\`ebre tensorielle sur $k$.\\

Les $({\cal{D}}_{X}^{(m)}) _{m \in {\Bbb{N}}}$ forment de mani\`ere naturelle un syst\`eme inductif dont la limite s'identifie \`a ${\cal{D}}_{X}$. Hormis lorsque $m=0$, l' anneau $ {\cal{D}}_{X}^{(m)} $ n'est pas engendr\'e par les d\'erivations mais admet localement (contrairement \`a ${\cal{D}}_{X}$) une famille finie de g\'en\'erateurs sur $O_{X}$.

\subsection{Alg\`ebres enveloppantes  de niveau $m$}  

Soit $G$ un $k$-groupe alg\'ebrique lin\'eaire d' alg\`ebre de Lie ${\cal{G}}$. Ici aussi, au lieu de d\'evelopper la th\'eorie des op\'erateurs diff\'erentiels de niveau $m \geq 0$  invariants sur $G$, il est plus rapide de  consid\'erer l'alg\`ebre ${\text{Dist}}\,(G)$ des distributions de $G$ (c' est \`a dire l'alg\`ebre des op\'erateurs diff\'erentiels invariants sur $G$)  et la filtration par l'ordre ${\text{Dist}}\,(G)= \bigcup_{n\in {\Bbb{N}}}{\text{Dist}}\,(G)_{n}$ et l' on pose pour $ m \geq 0$ ($\epsilon_{G}$ d\'esignant la co-unit\'e de $G$)\\

$U^{(m)}({\cal{G}})  : = {\Bbb{T}}_{k} ( {\text{Dist}}\,(G)_{2p^{m}-1}   ) / (\lambda-\lambda\epsilon_{G}, (\mu \otimes \mu '- \mu '\otimes \mu) -[\mu, \mu '] , \mu " \otimes \mu  -\mu " \mu  , \lambda \in k, \mu '' \in   {\text{Dist}}\,(G)_{ p^{m} }  ; \mu, \mu ' \in   {\text{Dist}}\,(G)_{ p^{m}-1}  )$\\

La limite inductive du syst\`eme inductif naturel que forment les  $(U^{(m)}({\cal{G}}))_{m\in {\Bbb{N}}})$ s'identifie \`a    ${\text{Dist}}\,(G)$ (qu' il est commode ici d' identifier \`a la ${\Bbb{Z}}$-forme de Kostant de l'alg\`ebre enveloppante de $sl_{2}$ r\'eduite modulo $p$ que nous noterons $U ({\cal{G}})$ dans la suite)  et $U^{(0)}({\cal{G}})$ n'est autre que l'alg\`ebre enveloppante universelle de ${\cal{G}}$.\\

Si d\'esormais $G = SL_{2}$ et si $(E,F,H)$ d\'esigne la base "standard"  de ${\cal{G}}$, la variante enti\`ere du th\'eor\`eme de Poincar\'e-Brikhoff-Witt implique que tout $U^{(m)}({\cal{G}})$ est engendr\'e sur $k$   par les \'el\'ements de la forme\footnote{Pour le lecteur non familier avec le yoga des puissances divis\'ees partielles, la formule (2.2.5.1) de \cite{Berthelot1} explique pourquoi les puissances divis\'ees partielles n'apparaissent pas ici bien qu'elles soient sous-jacentes.}\\

$E^{[a]}.(E^{[p^{m}]}) ^{a'}  \left(\begin{array}{c}H \\b\end{array}\right)  \left(\begin{array}{c}H \\p^{m}\end{array}\right) ^{b'} F^{[c]}.(F^{[p^{m}]})^{c'} $\\\\

avec $a,b,c \in [0, p^{m}-1]$, $a',b',c' \in {\Bbb{N}}$ (et $a!. E^{[a]} = E^{a}$, $b!. \left(\begin{array}{c}H \\b\end{array}\right):=  H(H-1)...(H-b+1 )$, $c!. F^{[c]} = F^{c}$). Quand $m \rightarrow + \infty$, on retrouve bien s\^ur la base "usuelle" sur $k$ de $U ({\cal{G}}) =  {\text{Dist}}\,(G) $.

\section{Scindage du Frobenius}

\subsection{La norme et ses propri\'et\'es}

La "norme" ${\Delta}_{T}$ (au sens de \cite{Haboush}, 6.) de l'alg\`ebre de Lie ${\cal{T}}$ du tore des matrices diagonales de $G$ est un \'el\'ement de $U({\cal{T}})$ qui peut se caract\'eriser   (\`a un facteur non nul pr\`es)  par des propri\'et\'es intrins\`eques (on renvoie \`a  {\it{loc. cit.}} pour les raisons de son introduction). Nous n'aurons   besoin que de sa description explicite et prendrons donc comme  d\'efinition (\cite{Haboush}, 6.3 Lemma) :\\

${\Delta}_{T} :=   \displaystyle\sum_{i=0}^{p-1} (-1)^{i} \left(\begin{array}{c}H \\i\end{array}\right) \in U ({\cal{T}}) \subset U ({\cal{G}})$ \\

Parmi les raisons de la terminologie "norme", on a \'evidemment  la propri\'et\'e, dans $U ({\cal{T}}) \subset U ({\cal{G}})$ :\\

  ${\Delta}_{T} ^{2} = {\Delta}_{T}$.\\

Pour la clart\'e des calculs \`a venir, il est \'egalement commode d'introduire, pour $l, i \in {\Bbb{Z}}$, les notations : \\

$\left(\begin{array}{c}H+l \\i\end{array}\right) := \frac{(H+l).(H+l-1)...(H+l-1)}{i!}$ \\

$\left(\begin{array}{c}   l   \\q\end{array}\right) =    \left(\begin{array}{c}H    \\q\end{array}\right)|_{ _{H=l}}$ (pour $l \geq 0$, c'est bien le coefficient binomial usuel) et\\

  ${\Delta} _{T,n} := \displaystyle\sum_{i=0}^{p-1} (-1)^{i} \left(\begin{array}{c}H-2n \\i\end{array}\right) \in U ({\cal{G}})$\\

 {\bf{Proposition 2.1.1}}. Soient $l  \in {\Bbb{Z}}$. On a, pour tout $r \in {\Bbb{N}}$, l'\'egalit\'e \\

 $  \left(\begin{array}{c}H+ l   \\r\end{array}\right) = \sum _{s+q=r, s\geq 0, q\geq 0}   \left(\begin{array}{c}   l   \\q\end{array}\right)   \left(\begin{array}{c}H    \\s\end{array}\right)  \in U ({\cal{G}})$.\\
 
 {\it{D\'emonstration}}. On d\'etermine les coefficients de $  \left(\begin{array}{c}H+ l   \\r\end{array}\right) $ dans la base des $(\left(\begin{array}{c}H    \\s\end{array}\right))_{s \geq 0}$ en faisant successivement $H=0, H=1,...$.\\

Cette \'egalit\'e est en fait valide dans ${\Bbb{Z}}[(\left(\begin{array}{c}H    \\s\end{array}\right))_{s \geq 0}]$ et posant   $H = H'+m$, on en d\'eduit  le\\
 
 {\bf{  Corollaire 2.1.2}}. Soient $l  \in {\Bbb{Z}}$ et $m \in{\Bbb{Z}}$. On a, pour tout $r \in {\Bbb{N}}$, l'\'egalit\'e \\

 $  \left(\begin{array}{c}H+ l +m  \\r\end{array}\right) = \sum _{s+q=r, s\geq 0, q\geq 0}   \left(\begin{array}{c}   l   \\q\end{array}\right)   \left(\begin{array}{c}H+m    \\s\end{array}\right)   $.\\

  {\bf{Corollaire 2.1.3}}\footnote{Je remercie X. Caruso de m'avoir signal\'e la formulation du corollaire suivant.}. On a ${\Delta}_{T} =   \left(\begin{array}{c}H- 1   \\p-1\end{array}\right) =   1 + \frac{H^{p-1}}{(p-1)!} = 1-H^{p-1}  $. En particulier, ${\Delta}_{T}(H)={\Delta}_{T}(aH)$ pour tout $a\in [0,...,p-1]$. Plus g\'en\'eralement, on a : ${\Delta} _{T,n}  =   \left(\begin{array}{c}H- 2n-1 \\p-1\end{array}\right) \in U ({\cal{G}})$\\
  
 {\it{D\'emonstration}}. Si $l = -1$ et $r= p-1$, la proposition 2.1.1 donne bien que ${\Delta}_{T} =   \left(\begin{array}{c}H- 1   \\p-1\end{array}\right) $. L'\'egalit\'e $ \left(\begin{array}{c}H- 1   \\p-1\end{array}\right) =   1 + \frac{H^{p-1}}{(p-1)!}       $   d\'ecoule du d\'eveloppement explicite du premier membre. Le reste est imm\'ediat.\\\\\
 
 {\bf{Corollaire 2.1.4}}.  On a, pour tout $n \in {\Bbb{Z}}$,  ${\Delta}_{T,n} ^{2} = {\Delta}_{T,n}$.\\

  {\bf{Corollaire 2.1.5}}.   L'\'el\'ement ${\Delta} _{T,n} \in U ({\cal{G}})$ ne d\'epend que de la classe modulo $p$ de $n$.\\

 {\bf{Proposition 2.1.6}}. On a, pour tout $n \geq 0$ et tout $m \in {\Bbb{Z}}$, les  \'egalit\'es \\ 
 
  $E^{[n]}.  \Delta_{T ,m} = \Delta_{T , m+n} . E^{[n]}$  \\
  
  $F^{[n]}.  \Delta_{T ,m} = \Delta_{T , m-n} . F^{[n]}$  \\
  
  {\it{D\'emonstration}}. cf. \cite{Haboush}, 6.4 Corollary ii).\\
  
  {\bf{Corollaire  2.1.7}}. On a :  $[E^{[np]}, \Delta_{T}] = [F^{[np]}, \Delta_{T}] = 0 $ pour tout $n \in {\Bbb{N}}$. \\

{\bf{Proposition 2.1.8}}. Si $p$ ne divise pas $j$, on a, pour tout $n \in{\Bbb{Z}}$ : $ \Delta_{T , n} .  \left(\begin{array}{c}H-2n    \\j\end{array}\right)  = 0$ sauf si $j$ est un multiple de $p$.


 {\it{D\'emonstration}}. On utilise que   ${\Delta} _{T,n}  =   \left(\begin{array}{c}H- 2n-1 \\p-1\end{array}\right)$ et l'on d\'etermine, comme pour la proposition 2.1.1,  les coordonn\'ees de $ \Delta_{T , n} .  \left(\begin{array}{c}H-2n    \\j\end{array}\right)$    dans la base $( \left(\begin{array}{c}H    \\i\end{array}\right))_{i\geq 0}$ en  faisant successivement $H=0$, $H=1$,.... .\\

\subsection{Le scindage}

Tout d'abord, rappelons que l' on dispose du morphisme de Frobenius $Fr : U ({\cal{G}}) \rightarrow U ({\cal{G}})$ (qui est le tranpos\'e  du morphisme de Frobenius usuel sur $G$) d\'efini sur les g\'en\'erateurs   par :\\

$Fr (E^{[n]}) = E^{[n/p]}$ , $Fr (F^{[n]}) = F^{[n/p]}$ \\


pour $n$ divisible par $p$ et $0$ sinon. C'est un morphisme d' alg\`ebres.\\

On va d\'efinir, pour tout $m \geq 0$ un  morphisme d'alg\`ebres    $\varphi_{m} : U^{(m )} ({\cal{G}}) \rightarrow U^{(m+1 )} ({\cal{G}})$  qui sera   induit par :\\

$\varphi (E^{[i]}) = E^{[ip]}.{\Delta}_{T}$ , $\varphi (F^{[i]}) = F^{[ip]}.{\Delta}_{T}$,  $\varphi ( \left(\begin{array}{c}H    \\i\end{array}\right)) =  \left(\begin{array}{c}H   \\ip\end{array}\right).{\Delta}_{T}$\\

avec $i \in {\Bbb{N}}$ (d'apr\`es le corollaire 2.1.7, on peut tout aussi bien multiplier \`a gauche ou \`a droite par ${\Delta}_{T}$). Cette application  scindera de mani\`ere \'evidente l'application $Fr$.  On a\\

{\bf{Proposition 2.2.1}}. L'application $\varphi_{m} : U^{(m )} ({\cal{G}}) \rightarrow U^{(m+1 )} ({\cal{G}})$          est bien d\'efinie et est un morphisme (non unif\`ere) d' alg\`ebres.\\

{\it{D\'emonstration}}. On va simplement d\'efinir un morphisme  $\varphi : U ({\cal{G}}) \rightarrow U ({\cal{G}})$ et laisser au lecteur le soin de v\'erifier qu' il se raffine en des morphismes $\varphi_{m}$ ; en effet,  formellement  tout \'el\'ement de $U^{(m )} ({\cal{G}}) $ est   multiple\footnote{cf. l'analogue de la formule (2.2.5.1) de \cite{Berthelot1}.} d'   \'el\'ements de $U  ({\cal{G}}) $ et si l'on sait d\'efinir $\varphi $ on sait d\'efinir $\varphi_{m}$. On \'evite ainsi le recours aux notations des puissances divis\'ees partielles qui ne jouent \`a cet endroit  aucun r\^ole. L'application $\varphi$ est la compos\'ee de l'application $Fr' : U ({\cal{G}}) \rightarrow U ({\cal{G}})$  ($Fr' (E^{[i]}) = E^{[ip]} $ , $Fr' (F^{[i]}) = F^{[ip]} $,  $Fr' ( \left(\begin{array}{c}H    \\i\end{array}\right)) =  \left(\begin{array}{c}H   \\ip\end{array}\right) $  apparaissant  dans la th\'eorie du "scindage de Frobenius" (voir \cite{Kumar2} 2. pour une pr\'esentation et \cite{Kumar1})  et de la multiplication par ${\Delta}_{T}$. Comme la restriction de l' application $Fr'$   \`a $U ({\cal{B}}^{+})$ ou bien \`a $U ({\cal{B}}^{-})$ (avec ${\cal{B}}^{+}$ resp. ${\cal{B}}^{-}$ l'alg\`ebre de Lie du sous-groupe des matrices triangulaires sup\'erieures, resp. triangulaires inf\'erieures de $G$) est un morphisme d'alg\`ebres, c'est {\it{a fortiori}} le cas apr\`es multiplication par  ${\Delta}_{T}$ car $ {\Delta}_{T}^{2} = {\Delta}_{T}$ et d' apr\`es le coroollaire 2.1.7. Il ne reste donc \`a v\'erifier que $\varphi$ (d\'efinie comme ci-dessus sur les g\'en\'erateurs de $U ({\cal{G}})$) est compatible avec la relation (cf. \cite{Humphreys}, 26.2 Lemma) :\\

$E^{[b]}F^{[a]} - F^{[a]}E^{[b]}= \sum_{r=1}^{{\text{min}}(a,b)} F^{[a-r]}  \left(\begin{array}{c}H-a-b+2r   \\r\end{array}\right) E^{[c-r]}$\\

On est donc amen\'e \`a consid\'erer :\\ 

$(E^{[pb]}F^{[pa]} - F^{[pa]}E^{[pb]})  .{\Delta}_{T}=   (\sum_{r=1}^{{\text{min}}(pa,pb)} F^{[pa-r]}  \left(\begin{array}{c}H-pa-pb+2r   \\r\end{array}\right) E^{[pb-r]}).{\Delta}_{T}  $\\

La proposition 2.1.6 donne $E^{[pb-r]} .{\Delta}_{T} =  \Delta_{T , pb-r} . E^{[pb-r]} $ et pour d\'emontrer la proposition,  il suffit donc de prouver que :\\

$\left(\begin{array}{c}H-pa-pb+2r   \\r\end{array}\right).\Delta_{T , pb-r} = \varphi(\left(\begin{array}{c}H- a- b+2r'   \\r'\end{array}\right))   $   si $r=r'p$  et \\\\

 $\left(\begin{array}{c}H-pa-pb+2r   \\r\end{array}\right).\Delta_{T , pb-r}  = 0$  si $r$ n' est pas divisible par $p$ \\

Prouvons la premi\`ere \'egalit\'e :\\

on a $\Delta_{T , pb-r} = \Delta_{T   } $ si $r=r'p$, d' apr\`es le corollaire 2.1.5. Il suffit donc de prouver que \\\\

$\left(\begin{array}{c}H-pa-pb+2r   \\r\end{array}\right).\Delta_{T } = \varphi(\left(\begin{array}{c}H- a- b+2r'   \\r'\end{array}\right)    $\\\\

 et donc que :\\


$Fr'(\left(\begin{array}{c}H- a- b+2r'   \\r'\end{array}\right) = \left(\begin{array}{c}H-pa-pb+2r   \\pr'\end{array}\right)$. \\

Pour calculer le membre de gauche,  on \'ecrit, gr\^ace \`a la proposition 2.1.1,$\left(\begin{array}{c}H- a- b+2r'   \\r'\end{array}\right) = \sum_{q+s=r', q\geq0, s\geq0} \left(\begin{array}{c}-a-b+2r'    \\q\end{array}\right) \left(\begin{array}{c}H    \\s\end{array}\right) $ dans la base des $\left(\begin{array}{c}H    \\j\end{array}\right) $, puis on applique $Fr'$.  Pour le membre de droite, on applique la m\^eme proposition et l' on remarque que les $\left(\begin{array}{c}m    \\j\end{array}\right) $ v\'erifient les  congruences usuelles :     $\left(\begin{array}{c}pm    \\j\end{array}\right)   \equiv 0 \,{\text{mod}}\,\, p$ sauf si $j$ est un multiple de $p$ et on a  $\left(\begin{array}{c}pm    \\j\end{array}\right) \equiv  \left(\begin{array}{c} m    \\j'\end{array}\right) \,{\text{mod}}\,\, p$ si $j=pj'$ (comme  cons\'equence imm\'ediate des "congruences   de  Lucas" : si $n=n_{d}p^{d}+...+n_{1}p+n_{0}$ avec $0  \leq   n_{i} < p$, $\left(\begin{array}{c}n    \\m\end{array}\right) \equiv  \left(\begin{array}{c}n_{d}    \\m_{d}\end{array}\right) .... \left(\begin{array}{c}n_{0}    \\m_{0}\end{array}\right) { \text{mod}} \,\, p$.) \\ 

Prouvons maintenant la deuxi\`eme  \'egalit\'e :\\

on a (coroll. 2.1.5) $\Delta_{T , pb-r}  = \Delta_{T , -r}$ et l' on a (Prop. 2.1.1)\\

 $ \left(\begin{array}{c}H-pa-pb+2r   \\r\end{array}\right)   = \sum _{s+q=r, s\geq 0, q\geq 0}   \left(\begin{array}{c}  -pa-pb   \\q\end{array}\right)   \left(\begin{array}{c}H+2r    \\s\end{array}\right)   $.\\
 
 Si l'on multiplie \`a droite par $\Delta_{T , -r}$, on d\'eduit de la proposition 2.1.8 que les termes de la comninaison lin\'eaire  que l' on obtient avec $s$ non divisible par $p$ sont nuls, d' autre part, si $s$ est divisible par $p$, $q$ ne peut l' \^etre (sinon $s+q=r$ le serait) et dans ce cas $\left(\begin{array}{c}  -pa-pb   \\q\end{array}\right)  \equiv 0 \,{\text{mod}}\, p$.\\
 
 {\bf{Remarque 2.2.2}}. Une   partie seulement (model\'ee sur \cite{Garnier}, Prop. 2.5.3 (3))   de ces calculs s'\'etend \`a des situations relev\'ees (c'est \`a dire modulo $p^{2}$,...) ou \`a l'alg\`ebre enveloppante quantique en une racine de l'unit\'e. \\

\section{Compatibilit\'es}

Avec les descriptions donn\'ees en 1., il est \'evident que l' application canonique $\rho  : U ({\cal{G}})    \rightarrow H^{0}(X,  {\cal{D}}_{X} )$ induit des applications $\rho_{m} : U^{(m)}({\cal{G}})    \rightarrow H^{0}(X,  {\cal{D}}_{X}^{(m)})$ et l'on va \'etablir la compatibilit\'e  entre celles-ci et les  deux applications $Fr  :  U^{(m+1)}({\cal{G}}) \rightarrow U^{(m)}({\cal{G}}) $et  $ \varphi_{m} : U^{(m)}({\cal{G}}) \rightarrow U^{(m+1)}({\cal{G}})$.\\

\subsection{Compatibilit\'e  avec l'application canonique de Berthelot}

L'application en question est celle d\'efinie dans  \cite {Berthelot3}, 2.1.7 ,  \cite {Berthelot2}, 2.2.4  : $can :  {\cal{D}}^{(m+1)}_{X} \rightarrow  F^{*}{\cal{D}}^{(m)}_{X' }$. Elle envoie, avec les notations de ${\it{loc. \,cit.}}$, $\partial^{<l>_{m+1}}$ sur $1 \otimes \partial '^{<l/p>_{m }}$ (avec la convention usuelle). On a alors la :\\

 {\bf{Proposition 3.1.1}}. Le diagramme suivant est commutatif

$$\xymatrix{
 U^{(m)}({\cal{G}}) \ar[d]_-{1\otimes \rho_{m }}  
 & U^{(m+1)}({\cal{G}}) \ar[d]^-{\rho_{m+1}} \ar@<0.5ex>[l]_-{Fr} \\
 H^{0}(X, F^{*}{\cal{D}}^{(m)}_{X' })  
 & H^{0}(X, {\cal{D}}^{(m+1)}_{X})  \ar@<0.5ex>[l]^-{can} }$$

 {\it{D\'emonstration}}. Pour v\'erifier cette commutativit\'e, il suffit  de   restreindre l' image par $\rho_{m+1}$ des \'el\'ements de $U^{(m+1)}({\cal{G}})$   \`a une  des deux copies de ${\Bbb{A}}^{1}$ formant le recouvrement standard de $X$ : $X = {\Bbb{A}}^{1} \cup {\Bbb{A}}^{1} = {\text{Spec}}\,(k[t]) \cup {\text{Spec}}\,(k[t'])$ car $H^{0}(X, {\cal{D}}^{(m)}_{X}) \hookrightarrow  H^{0}({\Bbb{A}}^{1}, {\cal{D}}^{(m)}_{{\Bbb{A}}^{1}})$. On peut donc supposer, avec les notations \'evidentes, que $E^{<n>_{m+1}}$ s' envoie sur $\partial_{t}^{<n>_{m+1}}$ (ou en passant \`a l' autre   ${\Bbb{A}}^{1}$ que $F^{<n>_{m+1}}$ s' envoie sur $\partial_{t'}^{<n>_{m+1}}$). La commutativit\'e est alors imm\'ediate. Comme ces \'el\'ements engendrent l'alg\`ebre $U^{(m)}({\cal{G}})$, cela suffit. \\
 
 {\bf{Remarque 3.1.2}}.  On remarquera que $\rho_{m+1} (\left(\begin{array}{c}H \\j\end{array}\right)_{m+1})$ est combinaison lin\'eaire de $t^{l}\partial^{<l>_{m+1}}$ et que lorsque $l=pl'$, on a $can(t^{l}\partial^{<l>_{m+1}}) = t^{l}(1 \otimes \partial '^{<l'>_{m }} ) = 1 \otimes t^{l'}\partial '^{<l'>_{m }}$.

\subsection{Compatibilit\'e  avec l'application d\'efinie par  Garnier}

D'autre part, dans \cite{Garnier}, 2.4, Garnier a introduit, pour $m > 0$,  un op\'erateur diff\'erentiel d' ordre $p-1$ \footnote{Cet op\'erateur est d\'efini de mani\`ere intrins\`eque, ie. sans aucun choix de coordonn\'ees, et donc "canonique".} (op\'erateur de Dwork) ${\cal{H}} \in {\cal{D}}^{(m)}_{X}$ pour tout sch\'ema lisse $X$ sur $k$ r\'ealisant la descente par Frobenius pour les ${\cal{D}}$-modules. Pour $X= {\Bbb{P}}_{k}^{1}\, (=X')$, cet op\'erateur s'\'etend  en fait une section globale\footnote{La question de savoir si tel est le cas pour toute vari\'et\'e de drapeaux semble int\'eressante.} de ${\cal{D}}^{(m)}_{X}$ (avec $m > 0$) sur $X$ et l' on a tout simplement :\\

{\bf{Proposition 3.2.1}}. On a $\rho_{m} (\Delta_{T}) =  {\cal{H}}$ dans $H^{0}(X, {\cal{D}}^{(m)}_{X})$ pour tout $m > 0$.\\

{\it{D\'emonstration}}. Soit $t$ une uniformisante locale sur $X$ au voisinage de $\infty$. On sait alors que ${\cal{H}} = \sum_{r=0}^{p-1}t^{r}\partial_{t}^{[r]}$. D'autre part, on a $\rho_{m} (H) = 2t\partial_{t}$ et l'on v\'erifie imm\'ediatement par r\'ecurrence que  $\left(\begin{array}{c}t\partial_{t} \\r\end{array}\right) = t^{r}\partial_{t}^{[r]}$. Le corollaire 2.1.3 permet alors de conclure.  \\

De plus, dans \cite{Garnier}, 4.3 et 4.4  est d\'ecrit un morphisme de $O_{X'}$-alg\`ebres ${\cal{D}}_{X'} \rightarrow F_{*}{\cal{D}}_{X}$    canonique (ie. ind\'ependant du choix de coordonn\'ees) lequel induit donc un morphisme $P \in  H^{0}(X, {\cal{D}}^{(m)}_{X}) \rightarrow P' \in  H^{0}(X, {\cal{D}}^{(m+1)}_{X})$ tel que $F'^{*} (P.f) = P'.(F'^{*}f)$ pour toute section $f \in O_{X'}$. Le lecteur prendra garde que, si $t$ d\'esigne une coordonn\'ee locale sur $X'$, contrairement \`a l'inclination naturelle,   $(1)' = {\cal{H}}\, (\neq 1)$ , $(\partial_{t})'  = \sum_{r=0}^{p-1} (-t)^{r} \left(\begin{array}{c}p+r \\p\end{array}\right)\partial^{[p+r]}_{t} \,(\neq \partial^{[p]}_{t})$. \\

 {\bf{Proposition 3.2.2}}. Le diagramme suivant est commutatif

$$\xymatrix{
 U^{(m)}({\cal{G}}) \ar[d]_-{\rho_{m}} \ar@<0.5ex>[r]^-{\varphi_{m}}
 & U^{(m+1)}({\cal{G}}) \ar[d]^-{\rho_{m+1}}  \\
 H^{0}(X, {\cal{D}}^{(m)}_{X}) \ar@<0.5ex>[r]_-{P \rightarrow P'}
 & H^{0}(X, {\cal{D}}^{(m+1)}_{X}) }$$

{\it{D\'emonstration}}. On effectue  la m\^eme r\'eduction \`a un calcul local que pour la proposition 3.1.1 et la proposition d\'ecoule alors imm\'ediatement des deux \'egalit\'es (la premi\`ere est un cas particulier de \cite{Garnier}, Prop. 4.6.2 (4), la seconde se v\'erifie par un calcul explicite) :\\

   $(\partial_{t}^{[n]})'  = \sum_{r=0}^{p-1} (-t)^{r} \left(\begin{array}{c}np+r \\np\end{array}\right)\partial^{[np+r]}_{t} $ \\
   
   $(\partial_{t}^{[np]}). {\cal{H}} = \sum_{r=0}^{p-1} (-t)^{r} \left(\begin{array}{c}np+r \\np\end{array}\right)\partial^{[np+r]}_{t} $\\

 \subsection{Descente des  $U^{(m)}({\cal{G}}) $-modules}
 
 J'ignore pour l'instant si l'on peut formuler un \'enonc\'e g\'en\'eral de descente pour ceux-ci. On peut toutefois essayer de s'inspirer de la Prop. 3.3.1 de \cite{Garnier} et, partant d' un $U^{(m+1)}({\cal{G}}) $-module $M_{m+1}$, consid\'erer ${\text{Im}}\,(M_{m+1}   \xrightarrow{m \rightarrow \Delta_{T} m} M_{m+1})$. \\
 
 Dans le cas du module de Verma infinit\'esimal $M_{m+1} := Z_{m+1}(-2)$ (dont la localisation, au sens de \cite{Bezrukavnikov}, correspond \`a la cohomologie \`a support en point \`a l'infini de $X$ du faisceau structural $O_{X}$), cette image  s'identifie exactement aux \'el\'ements sur lesquels $E,F,H$ agissent trivialement comme il r\'esulte de la description explicite de celui-ci (que l'on d\'eduit imm\'ediatement de   \cite {Mowbray}, 1.2). Une \'etude similaire \`a celle introduite dans ce travail devrait \'egalement tenir compte de la $p$-filtration sur  ${\cal{D}} _{X}$  et des th\'eories de localisation correspondantes (cf. \cite{Kaneda}). Nous reviendrons sur ces questions ult\'erieurement. \\
 
 {\bf{Remarque 3.3.1}}. En ce qui concerne le cas d' un groupe r\'eductif $G$ plus g\'en\'eral, remarquons que l' on a une injection $U^{(m)}({\cal{G}}) \hookrightarrow H^{0}(G, {\cal{D}}^{(m)}_{G})$ dont l' image consiste en les op\'erateurs diff\'erentiels invariants. Comme C. Noot-Huyghe me l' a fait remarquer, si l' on savait que les constructions de \cite{Garnier} sont $G$-\'equivariantes (par exemple au sens de \cite{Kaneda-Ye}), on   en d\'eduirait certainement    l' existence de $\varphi_{m}$ pour ces $G$. En ce qui concerne l' explicitation de ce dernier par les techniques de \cite{Garnier}, il est n\'eanmmoins \`a noter que si l' on veut travailler en \'ecrivant   $SL_{2} = {\text{Spec}}\,(k[a,b,c,d]/(ad-bc-1))$, la section $\Delta_{T}$ a une expression nettement plus compliqu\'ee (d\'eduite de $H \rightarrow a\partial_{a} - b\partial_{b} + c\partial_{c} - d\partial_{d}$) que les expressions apparaissant ci-dessus.

 \end{document}